\documentclass[oneside,notitlepage,12pt]{article}

\usepackage{amssymb}
\usepackage[leqno]{amsmath}
\usepackage{amsfonts}
\usepackage{amsopn}
\usepackage{amstext}
\usepackage{amsthm}
\usepackage{enumitem}
\usepackage[all]{xy}
\newdir{ >}{{}*!/-9pt/@{>}}

\usepackage[colorlinks, backref]{hyperref}
\usepackage{calrsfs}
\usepackage{fourier-orns}
\usepackage{hieroglf} 
\usepackage{clock} %\ClockFrametrue\ClockStyle2
\providecommand{\cal}{\mathcal}
\renewcommand{\Bbb}{\mathbb}

\newenvironment{pf}{\begin{proof}}{\end{proof}}

\newcommand{\Aaa}{{\cal{A}}}

\newcommand{\Cee}{{\cal{C}}}

\newcommand{\Ef}{{\cal{F}}}

\newcommand{\Pee}{{\cal{P}}}

\newcommand{\Yu}{{\cal{U}}}

\newcommand{\Nat}{{\Bbb{N}}}
\newcommand{\N}{{\Bbb{N}}}

\newcommand{\al}{\alpha}

\renewcommand{\phi}{\varphi}
\renewcommand{\rho}{\varrho}

\newcommand{\rest}{\restriction}

\newcommand{\ntr}{{n\in\omega}}

\newcommand{\loe}{\leqslant}
\newcommand{\goe}{\geqslant}

\newcommand{\subs}{\subseteq}
\newcommand{\sups}{\supseteq}

\newcommand{\setof}[2]{\{#1\colon #2\}}

\newcommand{\sett}[2]{\{#1\}_{#2}}
\newcommand{\sn}[1]{\{#1\}} % singleton
\newcommand{\dn}[2]{\{#1,#2\}} % doubleton
\newcommand{\pair}[2]{\langle #1, #2 \rangle} % pair
\newcommand{\map}[3]{#1\colon #2 \to #3} % A function
\newcommand{\power}[1]{\Pee(#1)}
\newcommand{\dpower}[2]{[#1]^{#2}}
\newcommand{\fin}[1]{[#1]^{<\omega}}

\newcommand{\bd}{{\mathfrak b}}

\newtheorem{tw}{Theorem}[section]
\newtheorem{wn}[tw]{Corollary}
\newtheorem{lm}[tw]{Lemma}
\newtheorem{prop}[tw]{Proposition}

\theoremstyle{definition}
\newtheorem{df}[tw]{Definition}
\newtheorem{ex}[tw]{Example}
\newtheorem{pyt}[tw]{Question}

\theoremstyle{remark}

\begin{document}

\title{On a topological Ramsey theorem}
\author{
{Wies{\l}aw Kubi\'s and Paul Szeptycki}
}
\date{\today\ \clocktime}

%\makeindex

\maketitle
\begin{abstract}
We introduce natural strengthenings of sequential compactness, the $r$-Ramsey property for each natural number $r\geq 1$.  We prove that metrizable compact spaces are $r$-Ramsey for all $r$ and give examples of compact spaces that are $r$-Ramsey but not $(r+1)$-Ramsey for each $r\geq 1$ (assuming CH for all $r>1$). Productivity of the $r$-Ramsey property is considered. 
\end{abstract}

\section{Introduction}

Let $K$ be a compact space and let $r$ be a positive integer.
Following \cite{BojKopTor2012}, we say that a function $\map f {\dpower Sr} K$ \emph{converges} to $p \in K$ if for every neighborhood $U$ of $p$ there is a finite set $F$ such that $f \dpower {S\setminus F}r \subs U$.
Once this happens for some $p$, we say that $f$ is \emph{convergent}.

This notion, for $r=2$ was introduced in \cite{BojKopTor2012} where the special case of our Theorem \ref{mainth} was stated and proved. Their main motivation was to obtain idempotents in compact semigroups $K$ as limits of certain functions $f:[\omega]^2\rightarrow K$. We show that more general notion of a space satisfying the $r$-Ramsey property (Definition \ref{maindef}) given below is a quite natural strengthening of sequential compactness and the main motivation of this paper is to prove some basic facts about this class of spaces and describe some examples showing that $r$-Ramsey can be strictly weaker than $(r+1)$-Ramsey. 

Our topological terminology is standard and basic definitions and notions can be found in $\cite{E}$. Set theoretic notation and terminology including some background on Ramsey's Theorem can be found in \cite{Ku}. And for a more detailed analysis of almost disjoint families, $\Psi$ spaces and the ideals FIN${}^n$, we refer the reader to \cite{Hr}.

\section{A sequential Ramsey theorem}

\begin{tw}\label{mainth}
Let $(K, \rho)$ be a compact metric space and let $\map f {\dpower \N r} K$ be an arbitrary function, where $r>0$.
Then there exists an infinite set $B \subs \N$ such that $f \rest {\dpower B r}$ converges to some element of $K$.
\end{tw}

\begin{pf}
For each $\ntr$ choose a finite cover $\Yu_n$ of $K$ consisting of balls of radius $2^{-n}$.
Then $\Yu_n$ induces a finite coloring of $\dpower \Nat r$.
Inductively, choose infinite sets $\Nat \sups A_0 \sups A_1 \sups \cdots$ such that
$A_n$ is monochromatic for the coloring induced by $\Yu_n$ (here we have used the classical Ramsey's theorem).
Let $B$ be any infinite set such that $B \setminus A_n$ is finite for every $\ntr$.
By compactness, $f \rest \dpower A r$ is convergent.
\end{pf}

The result above motivates the following definition:

\begin{df}\label{maindef}
Let $X$ be a topological space and let $r \in \Nat$ be positive.
We shall say that $X$ has the \emph{$r$-Ramsey property} (or $X$ is an $r$-Ramsey space) if for every function $\map f {\dpower \Nat r} X$ there exist $p \in X$ and an infinite set $B \subs \Nat $ such that $f \rest {\dpower B r}$ converges to $p$.
We shall say that $X$ has the \emph{Ramsey property} if it has the $r$-Ramsey property for every positive $r \in \Nat$. We will say that the set $B$  is a {\em convergent subsequence of $f$}. 
\end{df}

Note that the $1$-Ramsey property is just the sequential compactness.
Recall that a topological space $X$ is \emph{sequentially compact} if every sequence in $X$ has a convergent subsequence.

\begin{prop}
Every space with the $r$-Ramsey property has the $(r-1)$-Ramsey property, whenever $r>1$.
In particular, every space with the $r$-Ramsey property for some $r>0$ is sequentially compact.
\end{prop}

\begin{pf}
Assume $X$ has the $r$-Ramsey property and fix $\map g {\dpower \Nat {r-1}}X$.
Define $\map f {\dpower \Nat {r}}X$ by setting $f(s)=g(s \setminus \sn{\max s})$.
Let $A \in \dpower \Nat \omega$ be such that $f \rest {\dpower A r}$ converges to $p \in X$.

Fix a neighborhood $U$ of $p$.
There is $F \in \fin A$ such that $f(t) \in U$ whenever $t \subs A \setminus F$ and $|t|=r$.
Thus, if $s \in \dpower {A \setminus F}{r-1}$ then $g(s) = f(s \setminus \sn{\max s}) \in U$.
\end{pf}

It is useful to note that if $f:[S]^r\rightarrow X$ converges in a sequentially compact $X$, then $f$ has a somewhat nice canonical subsequence. 

\begin{df} A function $f:[S]^r\rightarrow X$ being {\em $r$-nice} is defined recursively. A sequence is $1$-nice if it is convergent. For $r>1$, $f:[S]^r\rightarrow X$ is $r$-nice if
\begin{enumerate}
\item $f$ is convergent to, say, $x$,
\item for every $s\in [S]^{r-1}$ the sequence $\{f(s\cup\{n\}):n\in S\setminus s\}$ is convergent to a point $x_s$, and
\item $g:[S]^{r-1}$ defined by $g(s)=x_s$ is $(r-1)$-nice. \end{enumerate}
\end{df}
Note that the fact that $g$, as defined in the definition, is convergent and converges to $x$ follows from the convergence of $f$ to $x$.

\begin{lm} If $X$ is sequentially compact, then any convergent $f:[S]^r\rightarrow X$ has an $r$-nice convergent subsequence. 
\end{lm}

\begin{proof} By induction on $r$. If $r=1$ then this follows since $X$ is sequentially compact. Fix $r>1$ and $f:[S]^r\rightarrow X$ convergent to some point $x\in X$. Enumerate $[S]^{r-1}$ as $\{s_k:k\in \omega\}$. Define $S_k$ recursively. $S_0\subseteq S$ is such that $\{f(s\cup\{n\}):n\in S_0\setminus s_0\}$ converges to $x_{s_0}$. And in general $S_{k+1}\subseteq S_k$ is chosen such that $\{f(s_k\cup\{n\}):n\in S_{k+1}\setminus s_k\}$ converges to some $x_{s_{k+1}}$.

Now take $T$ to be a pseudo-intersection of the $S_n$'s and note that $g:[T]^{r-1}\rightarrow X$ defined by $g(s)=x_s$ is defined and, since $f$ converges to $x$, so $g$ also converges to $x$. And so by our induction hypothesis, $g$ has an $(r-1)$-nice subsequence and so the lemma is proven.
\end{proof}
Note that it is easy to prove by induction that the closure of the image of an $r$-nice convergent $f$ is countable and so

\begin{wn}\label{countable} If $f:[\omega]^r\rightarrow X$ and $X$ has the $r$-Ramsey property, then $f$ has a convergent subsequence $f\upharpoonright[T]^r\rightarrow X$ such that the closure of $f^{\prime\prime}[T]^r$ is countable. 
\end{wn}

\section{Spaces with the Ramsey property}

It is clear that the class of all topological spaces with the Ramsey property is stable under closed subspaces and continuous images.
The same applies to the $r$-Ramsey property.
In order to see that there are arbitrarily large spaces with the Ramsey property, consider the $\Sigma$-product
$$\Sigma(\kappa) = \setof{x \in [0,1]^\kappa}{|\setof{\al}{x(\al)\ne0}| \loe \aleph_0}.$$
Note that the closure of every countable subset of $\Sigma(\kappa)$ is compact and metrizable, therefore $\Sigma(\kappa)$ has the Ramsey property.
More generally, all monolithic countably compact spaces have the Ramsey property (recall that a space $X$ is \emph{monolithic} if the closure of every countable subset of $X$ is second countable).

Recall that the \emph{unboundedness number} $\bd$ is the minimal cardinality of a family $\Ef \subs \omega^\omega$ which is unbounded with respect to the almost domination $<^*$, where $f <^* g$ if $f(n) < g(n)$ for all but finitely many $\ntr$.

\begin{tw}\label{smallcharacter}
Every sequentially compact space of character $<\bd$ has the Ramsey property.
\end{tw}

\begin{pf}
We use induction on $r$.
Sequential compactness implies that the theorem is true for $r=1$.
Suppose $r>1$ and the statement is true for $r-1$.
Fix $\map f {\dpower \N r} X$.

Given $a \in \N$, let $f_a(s) = f (s \cup \sn a)$, where $s \in \dpower {\N\setminus \sn a} {r-1}$.

We construct inductively a sequence $a_0 < a_1 < \cdots$ in $\N$, a sequence $\sett {p_n}{n>0} \subs X$, and a decreasing sequence of infinite sets $\N = A_0 \sups A_1 \sups \cdots$ such that
\begin{enumerate}
\item[(i)] $a_n \in A_n$,
\item[(ii)] $f_{a_n} \rest \dpower {A_{n+1}}{r-1}$ converges to $p_n$,
\end{enumerate}
Suppose $n>0$ and $a_i$, $A_i$, and $\sett{p_i}{i<n-1}$ have been constructed for $i<n$.
Using the theorem for $r-1$, we find $p_{n-1} \in X$ and an infinite set $A_n \subs A_{n-1}$ such that $f_{a_{n-1}} \rest \dpower {A_n}{r-1}$ converges to $p_{n-1}$.
We set $a_n = \min (A_n \setminus (a_{n-1}+1))$.

Now let $M \in \dpower \N \omega$ be such that $\sett {p_n}{n \in M}$ is convergent to some $p \in X$
(here we have used the sequential compactness of $X$).
Let
$$B = \setof{a_n}{n \in M}.$$
Re-enumerating $A_n$ and $p_n$, we may assume that $M = \Nat$.

Note that the set $B$ has the following property: Given a neighborhood $U$ of $p$, there is $m(U) \in \omega$ such that $\sett{p_n}{n \goe m(U)} \subs U$.
Consequently, for every $n \goe m(U)$ there exists $\phi_U(n) \in \omega$ such that
$f(s) \in U$ whenever $s \in \dpower B r$ is such that $\min s = a_n$ with $n \goe m(U)$ and $\min(s \setminus \sn{a_n}) \goe \phi_U(n)$ (the last fact follows from (ii), because $s \setminus \sn{a_n} \subs A_{n+1}$).
Define $\phi_U(n)$ arbitrarily for $n<m(U)$.

Let $\Yu(p)$ be a fixed base at $p$ such that $|\Yu(p)| < \bd$.
Then $\sett{\phi_U}{U \in \Yu(p)} \subs \omega^\omega$ has cardinality $<\bd$, therefore
we can find a strictly increasing function $\psi \in B^\omega$ such that $\phi_U <^* \psi$ for every $U \in \Yu(p)$.
Now let $C = \setof{\psi(n)}{\ntr}$.
We claim that $C$ is as required.

Fix $U \in \Yu(p)$ and fix $s \in \dpower C r$ such that $\min s = a_k$, where $k > m(U)$ and $\psi(n) > \phi_U(n)$ for every $n \goe k$.
Then $f(s) = f_{a_k}(s \setminus \sn{a_k}) \in U$.
This shows that $f \rest \dpower B r$ converges to $p$.
\end{pf}

\begin{wn}
Let $X$ be a sequentially compact space in which the closure of every countable set is first countable.
Then $X$ has the Ramsey property.
\end{wn}

\begin{wn}
Every countably compact linearly ordered space has the Ramsey property.
\end{wn}

Concerning products, it is known that the product of countably many sequentially compact spaces is sequentially compact. The same proof also applies to show
\begin{tw} The product of countably many $r$-Ramsey spaces is $r$-Ramsey. 
\end{tw} 
\begin{proof} Suppose we are given $r$-Ramsey spaces $X_i$ and  
$$f:[\omega]^r\rightarrow \prod_{i\in \omega} X_i$$
Recursively choose sets $B_0\supseteq B_1\supseteq ...$ so that for each $i$, 
$\pi_i\circ f:[B_i]^r\rightarrow X_i$ converges to $x_i$. Then it is straightforward to see that if $B$ is any pseudo-intersection of the family $\{B_i:i\in \omega\}$ then $f:[B]^r\rightarrow \prod_{i\in \omega} X_i$ converges to $\langle x(i):i<\omega\rangle$. 
\end{proof}

The splitting number, ${\mathfrak s}$ can be characterized as the minimal $\kappa$ such that $2^\kappa$ is not sequentially compact (see \cite{vD}). The analogous cardinal characteristic of the continuum that characterizes the $r$-Ramsey property in Cantor cubes is ${\mathfrak {par}}$ and was introduced by Blass in \cite{Bl}:
\begin{df} A set $A\subseteq\omega$ is {\em almost homogeneous} for a partition $f:[\omega]^n\rightarrow 2$ if there is a finite $F\subseteq A$ such that $f$ is constant on $[A\setminus F]^n$. ${\mathfrak {par}}_n$ denotes the smallest cardinal $\kappa$ such that there is a family of partitions of $[\omega]^n\rightarrow 2$ of size $\kappa$ such that no infinite set is almost homogeneous for all of them simultaneously.
\end{df}

First notice that ${\mathfrak {par}}_1$ is just the splitting number ${\mathfrak s}$. Also, note that if we consider partitions into some finite number of pieces $k$, instead of $2$ pieces, we obtain the same cardinal. And, moreover, for all $n\geq 2$ we have that ${\mathfrak {par}}_n={\mathfrak {par}}_2$, in fact:

\begin{tw}[\cite{Bl}] For each $n\geq 2$, ${\mathfrak {par}}_n=\min\{{\mathfrak {b, s}}\}$.
\end{tw}

Now we prove that the minimal $\kappa$ such that $2^\kappa$ is not $r$-Ramsey is precisely $\mathfrak{par}_2$.

\begin{tw} For each $\kappa<\mathfrak{par}_2$, $2^\kappa$ is $r$-Ramsey for all $r\in \omega$. And $2^{\mathfrak{par}_2}$ is not $2$-Ramsey.
\end{tw}
\begin{pf} Fix $\kappa<\mathfrak{par}_2$ and fix $f:[\omega]^r\rightarrow 2^\kappa$. For each $\alpha<\kappa$ let $f_\alpha=\pi_\alpha\cdot f$. By the definition of $\kappa<\mathfrak{par}_2=\mathfrak{par}_r$, there is $B\subseteq \omega$ and for each $\alpha$ an $i_\alpha$ such that for every $\alpha<\kappa$ there is a finite set $F$ such that $f_\alpha$ is constant with value $i_\alpha$ on $[B\setminus F]^r$. This just means that $f:[B]^r\rightarrow 2^\kappa$ converges to $(i_\alpha)_{\alpha\in \kappa}$ as required. 

To complete the proof, to see that $2^\kappa$ is not $2$-Ramsey for $\kappa=\mathfrak{par}_2$, fix a family $\{f_\alpha:\alpha<\kappa\}$ so that no $B\subseteq \omega$ is almost homogeneous for all functions $f_\alpha$. Taking the product function $f=\prod f_\alpha$ we have, as above, that for no $B$ can $f:[B]^2\rightarrow 2^\kappa$ converge. 
\end{pf}

\begin{wn} If $\mathfrak{b}<\mathfrak{s}$ then $2^{\mathfrak b}$ is sequentially compact but not $2$-Ramsey. 
\end{wn}

\section{Examples}

In this section we give some examples of spaces with the $k$-Ramsey Property that do not have the $(k+1)$-Ramsey Property. The first example of a $1$-Ramsey (i.e., sequentially compact) not $2$-Ramsey is in ZFC, but for larger $k$ we assume CH.  All of the examples are of the form $K(A)=\alpha(\Psi(\Aaa))$, the Alexandrov-Urysohn compactum formed by taking the one-point compactification of the $\Psi$-space determined by an almost disjoint family $\Aaa$ of infinite subsets of $\omega$. If $\Aaa$ is a maximal almost disjoint (mad) family, then $K(\Aaa)$ is not even $2$-Ramsey. 

\begin{ex}\label{madexample}
Let $\Aaa$ be a maximal almost disjoint family in $\dpower \Nat \omega$.
Then the Alexandrov-Urysohn compactum $K(\Aaa)$ is a compact and sequentially compact space failing the $2$-Ramsey property.
\end{ex}

\begin{proof} Recall that $K(\Aaa)$ is the Stone space of the Boolean subalgebra of $\power \Nat$ generated by $\fin \Nat \cup \Aaa$. 
For convenience, we assume that $\Aaa \subs \dpower{\Nat \times \Nat}\omega$ and $\sn n \times \Nat \in \Aaa$ for each $n \in \Nat$.
Thus
$$K(\Aaa) = (\Nat \times \Nat) \cup \setof{p_A}{A \in \Aaa} \cup \sn \infty,$$
where we write $p_n$ instead of $p_{\sn n \times \Nat}$.
Then $\Nat \times \Nat$ is the isolated points, $K(\Aaa) \setminus \sn \infty$ is locally compact and a basic neighborhood of $p_A$ is $\sn {p_A} \cup (A \setminus F)$, where $F \in \fin A$.

Let $\map f{\dpower \Nat 2}{K(\Aaa)}$ be defined by $f(\dn k \ell) = \pair k \ell$, where $k < \ell$.
Now it is easy to see that this $f$ witnesses the failure of $2$-Ramsey since for any $S\subseteq \omega$ infinite, the closure of $f^{\prime\prime}[B]^2$ in $\Psi(\Aaa)$ is infinite, and since $\Aaa$ is mad, it follows that the closure is uncountable. Therefore by Corollary \ref{countable}, $K(\Aaa)$ can't be $2$-Ramsey.
\end{proof}
%We claim that there is no $B \in \dpower \Nat \omega$ such that $f \rest \dpower B2$ is convergent.
%Fix $B \in \dpower \Nat \omega$ and suppose that $f \rest \dpower B2$ converges to $p \in K(\Aaa)$. 

%We first note that $p\not=\infty$. Indeed, if $F$ is any finite set, then $\{(k,l)\in (B\setminus F)^2:k<l\}$ is in $I^+(\Aaa)$ and it is well known (see, e.g., see Proposition 2 in \cite{Hr}) that there is therefor a $A\in \Aaa$ such that $A\cap \{(k,l)\in (B\setminus F)^2:k<l\}$ is infinite for all finite $F\subseteq B$. I.e., $p$ can't be the point $\infty$ and so $p = p_A$ for some $A \in \Aaa$ and $A \ne \sn n \times \Nat$ for every $n \in \Nat$.
%Consider the neighborhood $U = \sn {p_A} \cup A$ of $p_A$ and fix a finite set $F \subs B$.
%Fix $k \in B \setminus F$.
%Then the infinite set
%$$\setof{\pair k n }{n\in B,\; n > k} \subs \sn n \times \Nat$$
%is almost disjoint from $A$, therefore we may find $\ell \in B$, $\ell > k$ such that $\pair k \ell \notin A$.
%Thus $f(\dn k \ell) \notin U$, a contradiction. 

%Note that if $A$ is mad, as in the above example, then $K(A)$ is not Fr\'echet. Indeed, $K(A)$ is Fr\'echet precisely when $A$ is nowhere mad \cite{Hr}, meaning that for every $X\subseteq \omega$, if $X$ is not  covered by finitely many sets from $A$, then there is a $b\subseteq X$ such that $b$ is almost disjoint from each $a\in A$. So to find a $K(A)$ that is Fr\'echet but not $2$-Ramsey, 

%Another example of an Alexandrof-Urysohn compactum that is not $2$-Ramsey 

We now describe, for each $r>1$ an $r$-Ramsey not $(r+1)$-Ramsey compact space of the form $K(\Aaa)$. We state a few lemmas about these properties in these types of spaces. 

\begin{lm}\label{lemma1} For any almost disjoint family $\Cee$ on a countable set $I$ of isolated points, $K(\Cee)$ is $r$-Ramsey if and only if for any $f:[\omega]^r\rightarrow I$, there is $B\subseteq \omega$ such that $f\upharpoonright [B]^r$ converges. 
\end{lm}

\begin{proof} The property is clearly necessary. For sufficiency, note that for any $f:[\omega]^r\rightarrow K(\Cee)$ one can find $B$ such that either $f''[B]^r\subseteq \Cee$ (in which case, since any infinite subset converges, one can easily find $B'\subseteq B$ witnessing $r$-Ramsey) or $f''[B]^r\subseteq I$, as required.
\end{proof}

\begin{lm}\label{lemma2} $K(\Aaa)$ is $r$-Ramsey if for every $f:[\omega]^r\rightarrow I$ there is $B\subseteq \omega$ such that 
$$
\Cee \upharpoonright f''[B]^r=\{p\in \Cee: p\cap f''[B]^r\text{ is infinite}\}
$$ 
is of size less than ${\mathfrak b}$.
\end{lm}
\begin{proof} If we can find such a $B$, then the closure of the subset $f''[B]^r$ in $K(\Aaa)$ has countable character at all points of $\Aaa$ and character less than ${\mathfrak b}$ at $\infty$ and so by Theorem \ref{smallcharacter} this subspace is $r$-Ramsey and so we can find $B'\subseteq B$ on which $f$ converges. 
\end{proof}

We now turn to a construction of an $r$-Ramsey example of an Alexandrov-Urysohn compactum that is not $(r+1)$-Ramsey. 

We first need some basic properties of the Fubini product of the ideal FIN. Recall, FIN is the ideal of finite subsets of $\omega$ and for each $n>1$, FIN${}^n$ is the ideal on $\omega^n$ defined recursively by 
$$A\in \text{FIN}^n \text{ if } \{\overline{k}\in \omega^{n-1}:\{j:\overline{k}\frown n\in A\} \text{ is infinite}\}\in \text{FIN}^{n-1}.$$

The following lemma is easily proved using the definition. 

\begin{lm} For any $X\subseteq \omega^n$, if $X\not\in\text{FIN}^n$ then there is a $T\subseteq X$ such that $T$ forms an everywhere $\omega$-splitting tree. I.e., letting $T_i=\{x\upharpoonright i: x\in T\}$, for each $i<n$ and each $s\in T_i$,  
$\{t\in T_{i+1}:t\upharpoonright i = s\}$ is infinite. 
\end{lm}

Now we need the following
\begin{lm}\label{lemmaFIN} For any function $f:[\omega]^n\rightarrow \omega^{n+1}$ there is a $B\subseteq \omega$ infinite such that $f''[B]^n\in \text{FIN}^{n+1}$.

\end{lm}
\begin{proof} By induction on $n$. When $n=1$ we can find $B$ such that $f''B$ is either contained in a column $\{k\}\times\omega$ or is a partial function. In either case $f''B\in {\text FIN}^2$. 

For the inductive step, fix $f:[\omega]^n\rightarrow \omega^{n+1}$ and fix $k_0\in \omega$ arbitrary. Define $f_{k_0}:[\omega\setminus k_0+1]^{n-1}\rightarrow \omega^{n+1}$ by $f_{k_0}(s)=f(\{k_0\}\cup s)$. By our inductive assumption, there is a $B_1\subseteq \omega\setminus k_0+1$ such that 
\begin{enumerate}
\item[$(*)$] the projection of $f_{k_0}''[B_1]^{n-1}$ onto any $n$ coordinates of $\omega^{n+1}$ is not in FIN${}^n$. 
\end{enumerate}
Let $k_1=\min(B_1)$ and continue recursively constructing $\{k_i,B_i\}$ so that for each $i$, $k_i=\min(B_i)$ and 
\begin{enumerate}
\item[$(*)$] the projection of $f_{k_i}''[B_{i+1}]^{n-1}$ onto any $n$ coordinates of $\omega^{n+1}$ is not in FIN${}^n$. 
\end{enumerate}
Now let $B=\{k_i:i\in \omega\}$. Since $(\omega+1)^{n+1}$ is $r$-Ramsey for all $r$, we may, by possibly shrinking $B$, assume that $f:[B]^n\rightarrow (\omega+1)^{n+1}$ converges to some $x$. 

CASE 1. $x(i)=\omega$ for all $i<n+1$. In this case, if $f''[B]^n$ is not in FIN${}^{n=1}$ then there is an everywhere $\omega$-splitting tree $T\subseteq f''[B]^n$ witnessing this. Fix $z\in T$ and consider $k=z(0)$. Since $f:[B]^n\rightarrow (\omega+1)^{n+1}$ converges to $x$, there is an $N$ such that 
$$f''[B\setminus N]^n\subseteq (\omega\setminus k+1)^{n+1}$$
 Therefore, $f''[B\setminus N]^n$ is disjoint from $T'=\{s\in T: s(0)<N\}$. Note that the projection of $T'$ onto $\omega^{n+1\setminus \{0\}}$ is not in FIN${}^n$. And $T'$ is the union of the sets of the form $f''\{s\in [B]^n:\min s < N\}$. But for each $k_i<N$ the projection of the sets $f''\{s\in [B]^n:\min s =k_i\}$ onto $\omega^{n+1\setminus\{0\}}$ is in FIN${}^n$. Contradiction.
 
CASE 2. There is $i$ such that $x(i)\not=\omega$. In this case, there is some $N$ such that 
$$
f''[B\setminus N]^n\subseteq \{z\in (\omega)^{n+1}: z(i)=x(i)\}.
$$
And this latter set is not in FIN${}^{n+1}$.
\end{proof}

%Assuming the lemma holds for $n-1$, fix $f:[\omega]^n\rightarrow \omega^{n+1}$. Let $k_0=0$ and fix $B_0\subseteq \omega$ such that  By a standard diagonalization argument, we can find $B\subseteq \omega$ infinite, such that for each $\overline{k}\in \omega^n$, $f''[B]^n \cap \overline{k}\times \omega$
 
Now fix $r>1$ and we will need to assume CH. We build an almost disjoint family $\Aaa$ on $\omega^{r+1}$ so that for the function $G:[\omega]^{r+1}\rightarrow \omega^{r+1}$ defined by $G(x)=\langle x(0),...,x(r)\rangle$ where $x=\{x(0),...,x(r)\}$ is the increasing enumeration of $x$, $G$ will have no convergent subsequence $B$. 

To simplify some notation for any $B$, let 
$$
B^{\uparrow n}=\{\overline{k}\in B^n: \overline{k}(i)<\overline{k}(i+1)\text{ for all }i<n-1\}.
$$

\begin{lm}\label{lemmanotramsey} Suppose $\Aaa$ is an almost disjoint family on $\omega^{r+1}$. For the function $G$ defined as above, $B$ is a convergent subsequence for $G$ with limit $\infty$ in $K(\Aaa)$ if and only if for every $a\in \Aaa$ there is an $n$ such that $a\cap (B\setminus n)^{\uparrow r+1}=\emptyset$. 
\end{lm}
\begin{proof} Directly from the definitions.
\end{proof}

\begin{tw} Assume CH. For each $r>1$ there is an almost disjoint family $\Aaa$ on $\omega$ such that $K(\Aaa)$ is $r$-Ramsey but not $(r+1)$-Ramsey. 
\end{tw}
\begin{proof} We will construct $\Aaa=\{a_\alpha:\alpha<\omega_1\}$ an almost disjoint family on $\omega^{r+1}$ by defining each $a_\alpha$ by recursion on $\alpha$.
We start by letting $\{a_n:n\in\omega\}$ be an enumeration of the disjoint family 
$$
\{\{\overline{k}\}\times\omega:\overline{k}\in \omega^r\}
$$
Note that each $a_n\in \text{FIN}^{n+1}$ and any $x$ that is almost disjoint from all the $a_n$ is also in FIN${}^{n+1}$. 

We enumerate as $\{B_\alpha:\alpha<\omega_1\}$ all infinite subsets of $\omega$ and fix an enumeration
$$
\{f_\alpha:\alpha<\omega_1\}=\left(\omega^{r+1}\right)^{[\omega]^r}.
$$
Recall our plan that $G:[\omega]^{r+1}\rightarrow \omega^{r+1}$ should have no convergent subsequence. 

Suppose that at some stage $\alpha$ of the construction we have defined
$\{a_\beta:\beta<\alpha\}$ and 
$\{X_\beta:\beta<\alpha\}$
such that
\begin{enumerate}
\item For all $\beta<\alpha$, $f_\beta''[X_\beta]^r\in \text{FIN}^{n+1}$.
\item For all $\beta<\gamma<\alpha$, $a_\beta\cap a_\gamma$ is finite. 
\item For all $\beta<\gamma<\alpha$, $\left(f_\beta''[X_\beta]^r\right)\cap a_\gamma$ is finite. 
\item For all $\beta<\alpha$ and all $n<\omega$, $a_\beta\subseteq^*G''[B_\beta\setminus n]^{r+1}$.
\end{enumerate}
To define $a_\alpha$ and $X_\alpha$, note first that the following family of sets
$${\mathcal H}=\{f_\beta''[X_\beta]^r:\beta<\alpha\}\cup\{a_\beta:\beta<\alpha\}\cup\{G''[B_\alpha\setminus n]:n\in \omega\}\subseteq \text{FIN}^{n+1}.$$
Therefore the family $\{(\omega^{n+1})\setminus H: H\in {\mathcal H}\}$ has an infinite pseudo-intersection. Let $a_\alpha$ be any such pseudo-intersection and it directly follows that $a_\alpha$ satisfies the inductive hypotheses (2)--(4). So we need only define $X_\alpha$ to satisfy (1). But that there is such an $X_\alpha$ follows from Lemma \ref{lemmaFIN}.

This completes the construction of $\Aaa=\{a_\alpha:\alpha<\omega_1\}$.
To see that it is $r$-Ramsey, by Lemma \ref{lemma1} we need only consider functions $f:[\omega]^r\rightarrow \omega^{r+1}$, and each such $f$ appears as an $f_\alpha$. For each $\alpha$ by inductive hypothesis (3) we have that 
$$
\{a\in \Aaa: a\cap f_\alpha''[X_\alpha]^r\text{ is infinite}\}
$$ 
is countable. And so by Lemma \ref{lemma2}, it follows that $K(\Aaa)$ is $r$-Ramsey. 
On the other hand, to see that $G$ has no convergent subsequence, we note that by inductive hypothesis (4) we have, by Lemma \ref{lemmanotramsey}, that no $B_\alpha$ is a convergent subsequence for $G$ and so $K(\Aaa)$ is not $(r+1)$-Ramsey.
\end{proof}

\section{Questions}

We finish with the following questions:

\begin{pyt}
Does there exist a space with the $r$-Ramsey property and without the $(r+1)$-property assuming only ZFC? Here $r>1$. 
\end{pyt}
Note that all our examples of $r$-Ramsey and not $(r+1)$-Ramsey spaces (which required CH) are Fr\'echet-Urysohn but the ZFC example that was not $2$-Ramsey, being a $K(\Aaa)$ where $\Aaa$ is mad, is not Fr\'echet. However, a similar example could be constructed from a completely separable mad family. The main idea is to start with a completely separable mad family on $\omega\times\omega$ as in the construction of Example \ref{madexample}. Then, deleting an infinite set from each $A\in\Aaa$ other than the fixed columns $\{n\}\times\omega$ will still give an example that fails to be $2$-Ramsey. But it will be Fr\'echet since it is nowhere mad ((see \cite{Hr} for the definitions of completely separable and nowhere mad). So we have 
\begin{ex} Assuming the existence of a completely separable mad family, there is an $\Aaa$ such that $K(\Aaa)$ is Fr\'echet and not $2$-Ramsey. 
\end{ex}
Although the existence of a completely separable mad family is a relatively weak one (e.g., it follows from ${\mathfrak c}<\aleph_\omega$ or ${\mathfrak s}\leq{\mathfrak a}$), we ask
\begin{pyt}
Does there exist a ZFC example of a Fr\'echet-Urysohn compact space without the Ramsey property?
\end{pyt}

We know that a product of any countable family of $r$-Ramsey spaces is $r$-Ramsey, and we have characterized when $2^\kappa$ is $r$-Ramsey. Moreover, ${\mathfrak h}$ is the minimal $\kappa$ such that a product of fewer than $\kappa$ many sequentially compact spaces is sequentially compact \cite{Si}, and we conjecture the same holds for $r$-Ramsey. 

Indeed, the proof that the product of fewer than ${\mathfrak h}$ sequentially compact spaces is sequentially compact also shows the same for $r$-Ramsey, but the family of sequentially compact spaces whose product is not sequentially compact given in \cite{Si} are, in fact, Alexandrov-Urysohn compacta that we have seen are not even $2$-Ramsey. So, if $\mu$ is the minimal cardinal for productivity of the class of $2$-Ramsey spaces, then $\mathfrak{h}\leq \mu$ and we conjecture that $\mu=\mathfrak{h}$. 

\begin{pyt} Characterize the minimal cardinal $\kappa$ satisfying the product of fewer than $\kappa$ many $r$-Ramsey spaces is always $r$-Ramsey.
\end{pyt}

%The collection of compact spaces with no non-trivial convergent sequences is a very interesting class of spaces and the famous Efimov problem asks whether it is consistent that every compact space either contains a non-trivial convergent sequence or a copy of $\beta\omega$. One might consider a convergent sequence of the form $f:[B]^r\rightarrow X$ r-nontrivial if $f$ is 1-1 or finite-1 and consider the class of spaces with no $r$-nontrivial convergent sequences and ask if they may have a nontrivial convergent sequence or an $r-1$-nontrivial convergent sequence. And is there an interesting and non-trivial formulation of Efimov's problem for $r$-nontrivial convergence? 

\paragraph{Acknowledgments.} The first author would like to thank Henryk Michalewski for bringing this topic and in particular for pointing out reference~\cite{BojKopTor2012}.\\
The second author would like to thank the Czech Academy of Sciences and NSERC for support.

\end{document}